\documentclass[a4paper,12pt]{amsart}
\usepackage{amssymb}
\usepackage[alwaysadjust]{paralist}
\usepackage[capitalize]{cleveref}
\usepackage{graphicx}
\usepackage{float}
\numberwithin{equation}{section}

\newtheorem{thm}[equation]{Theorem}
\theoremstyle{definition}

\newtheorem{rem}[equation]{Remark}
\newtheorem{rems}[equation]{Remarks}
\newtheorem{que}[equation]{Question}
\newcommand{\tr}{\operatorname{tr}}

\begin{document}
\title[angle-spectrum]{An arithmetic property of the set of angles between closed geodesics on hyperbolic surfaces of finite type}
\author{Sugata Mondal}
\address{Indiana University,
	Rawles Hall, 831 E 3rd Street,
	Bloomington, Indiana\\}
\email{sumondal\@@iu.edu}
\thanks{\emph{Acknowledgments.} I would like to thank Chris Judge for the discussions that I had with him on this problem. The example in the last section is due to Hugo Parlier that was explained to the author by Chris Judge.}
\date{\today}

\maketitle

\begin{abstract}
For a hyperbolic surface $S$ of finite type we consider the set $A(S)$ of angles between closed geodesics on $S$. Our main result is that there are only finitely many rational multiples of $\pi$ in $A(S)$. 
%A complete description of this finite set for the modular surface
%$\mathbb{H}^2/\text{PSL}_2(\mathbb{Z})$ is discussed.
\end{abstract}

\section{Introduction}
Geometry of two dimensional manifolds, surfaces, have been in the center of mathematical research for 
centuries. Hyperbolic metrics on surfaces has been a very important testing ground for different 
geometric curiosities. In particular, lengths of closed geodesics on these surfaces have been an important topic of research for years (see \cite{Hu}, \cite{O}, \cite{W}). This article is on a related geometric quantity the angles between pairs of closed geodesics.  

Let $S$ be a hyperbolic surface of finite type. We denote the set of angles between pairs of closed geodesics on $S$ by 
$\mathcal{A}(S)$. A fixed angle may appear at many different intersections. We call this number of distinct appearances the {\it multiplicity} of the angle. We denote the set of angles in $\mathcal{A}(S)$ forgetting their multiplicities by $A(S)$ and call $\mathcal{A}(S)$ and $A(S)$ by {\it angle spectrum} and {\it angle set} respectively. 

We begin by specifying a way of measuring these angles. Let $\gamma$ and $\delta$ be two closed geodesics on $S$ that intersect each other at $p$. We measure the angle of intersection $\theta(\gamma, \delta, p)$ in the counter clockwise direction from $\gamma$ to $\delta$. In particular
$\theta(\gamma, \delta, p) = \pi - \theta(\delta, \gamma, p)$.

\begin{rem}
For $\gamma, \delta$ and $\theta = \theta(\gamma, \delta, p)$ as above $\cos^2(\theta)$ depends on $\gamma$ and $\delta$ but not on the direction in which the angle
is measured. 	
\end{rem} 
In this article we focus on qualitative properties of the two collections $\mathcal{A}(S)$ and $A(S)$. For any hyperbolic surface $S$ of finite type $A(S)$ is a countable infinite set and it follows from \cite{P-S} that $A(S)$ is dense in $[0, \pi]$. The main question that we address in this paper is the following. 
\begin{que}\label{multi}
 How many angles in $\mathcal{A}(S)$ can be a rational multiple of $\pi$ ?
\end{que}
Surprisingly the author's motivation to study this question came from a seemingly unrelated field. In the paper \cite{J-M} we have studied eigenfunctions of the Laplacian on hyperbolic surfaces. Let $\phi$ be an eigenfunction of the Laplacian on $\mathbb{H}^2$. Let $\Gamma_\phi$ denote the subgroup of isometries of $\mathbb{H}^2$ that leaves $\phi^2$ invariant. In \cite{J-M} we have observed (motivated by a similar observation in \cite{G-R-S}) that if $\phi$ vanishes on a geodesic $\gamma$ then it is odd with respect to the reflection isometry $R_\gamma$ along $\gamma$ of $\mathbb{H}^2$ i.e. $\phi \circ R_\gamma = - \phi$. In particular, $R_\gamma \in \Gamma_\phi$. An important property of any non-constant eigenfunction $\phi$ is that the subgroup $\Gamma_\phi$ is discrete (see \cite{J-M}).

Now consider two intersecting geodesics $\gamma, \delta$ on $\mathbb{H}^2$ and consider the subgroup $\Gamma(\gamma, \delta)$ of SL$(2, \mathbb{R})$ generated by the reflections $R_\gamma, R_\delta$ along $\gamma$ and $\delta$ respectively. Let $\gamma, \delta$ intersect each other at $p$ and let $\theta = \theta(\gamma, \delta, p)$. Then $\Gamma(\gamma, \delta)$ contains an elliptic isometry of $\mathbb{H}^2$ which is a rotation about $p$ by an angle equal to $\theta$. Let $\phi$ be a non-constant eigenfunction that vanish on both $\gamma$ and $\delta$. Then by the last paragraph $\Gamma(\gamma, \delta) \subset \Gamma_\phi$. In particular, since $\Gamma_\phi$ is discrete, so is $\Gamma(\gamma, \delta)$ implying that $\theta$ must be a rational multiple of $\pi$. 

It is not difficult to construct eigenfunctions that vanish on two intersecting geodesics, even on closed hyperbolic surfaces (see \cite{J-M}). In general the answer to Question \ref{multi} is `infinite'. In the last section we construct examples of surfaces for which there are infinitely many distinct intersections between pairs of closed geodesics such that the angle of intersection is $\pi/2$. The main result of this article is that, in general, there are infinitely many rational multiples of $\pi$ in $\mathcal{A}(S)$ if and only if one of these rational multiples of $\pi$ has infinite multiplicity in $\mathcal{A}(S)$. 
\begin{thm}\label{fuchsian}
For any hyperbolic surface $S$ of finite type there are only finitely many angles in $A(S)$ that are rational multiples of $\pi$.
\end{thm}
\subsection{Structure of the article}
In the first section we consider a complete hyperbolic surface $S$ of finite type. Using uniformization theorem we consider a Fuchsian group $\Gamma$ such that $S = \mathbb{H}^2/\Gamma$, up to isometry. For two intersecting closed geodesics $\gamma$ and $\delta$ on $S$ we fix an intersection point $p$. In \S1 we consider $M_\gamma, M_\delta \in \Gamma$ representing $\gamma$ and $\delta$ respectively and use the matrix entries of $M_\gamma$ and $M_\delta$ to get a formula for $\cos^2 (\theta)$ where $\theta = \theta(\gamma, \delta, p)$. 

We prove Theorem \ref{fuchsian} in \S2. In the first step of the proof we consider the field $\mathbb{F}_\Gamma$ obtained by attaching the matrix entries of all the elements in $\Gamma$ to $\mathbb{Q}$. Using the fact that $\Gamma$ is finitely generated it follows that $\mathbb{F}_\Gamma$ is a finitely generated field extension of $\mathbb{Q}$. Using the expression for $\cos^2(\theta)$ obtained in \S1 we deduce that $\cos^2(\theta) \in \mathbb{F}_\Gamma$ for any angle $\theta \in A(\mathbb{H}^2 /\Gamma)$. 

The final arguments of the proof go follows. For simplicity, assume that $\mathbb{F}_\Gamma$ is algebraic over $\mathbb{Q}$. Since $\mathbb{F}_\Gamma$ is finitely generated over $\mathbb{Q}$ we obtain that the degree of extension $\mathbb{F}_\Gamma|_{\mathbb{Q}}$ is finite. Now let $\frac{p}{q}\pi$ be in $A(\mathbb{H}^2 /\Gamma)$ and so $\cos^2(\frac{p}{q}\pi) \in \mathbb{F}_\Gamma$. Then there is a field extension $\mathbb{F}(q)$ of $\mathbb{F}_\Gamma$ with degree of extension at most two that contain a primitive $q$-th root of unity. In particular, the degree of extension $\mathbb{F}(q)|_\mathbb{Q}$ is uniformly bounded independent of $q$. Finally we observe that the degree of extension $\mathbb{F}(q)|_\mathbb{Q}$ is at least $\phi(q)$ where $\phi$ is the Euler's $\phi$-function that counts the number of distinct positive integers less than and co-prime with $q$. Since $\phi(q)$ goes to infinity as $q$ goes to infinity \cite[Theorem 328]{H-W}, we reach our desired contradiction.
\section{Formula for the cosine of an angle}
Let $\Gamma$ be a finitely generated Fuchsian group. This usually means that $\Gamma \subset \textrm{PSL}(2, \mathbb{R})$. By taking the pre-image of $\Gamma$ under the quotient map $\Pi: \textrm{SL}(2, \mathbb{R}) \to \textrm{PSL}(2, \mathbb{R})$ we can always think of $\Gamma \subset  \textrm{SL}(2, \mathbb{R})$. This identification will be assumed in the article from now on. It is a standard fact that every closed geodesic on $S$ corresponds to a conjugacy class of elements in $\pi_1(S) = \Gamma$. Let $\gamma, \delta$ be two closed geodesics on $S$ and let $M_\gamma, M_\delta \in \Gamma$ be two representatives of $\gamma, \delta$ respectively. Let us denote
\begin{equation*}
 M_\gamma =   \left( \begin{array}{cc}
 a_\gamma & b_\gamma \\
 c_\gamma & d_\gamma \
 \end{array} \right),   ~~~~~                     M_\delta =   \left( \begin{array}{cc}
                                                  a_\delta & b_\delta \\
                                                  c_\delta & d_\delta \
                                                  \end{array} \right).
\end{equation*}
Recall that $\gamma$ and $\delta$ are the projections of the axes of  $M_\gamma$ and $M_\delta$ respectively, under the covering map: $\mathbb{H}^2 \to \mathbb{H}^2 /\Gamma$. Since  $\gamma$ and $\delta$ are closed geodesics, $M_\gamma$ and $M_\delta$ are hyperbolic linear fractional transformations. Thus the axes of $M_\gamma$ and $M_\delta$ are either semi-circles or vertical straight lines that intersect  $\mathbb{R}$ orthogonally. Here $\mathbb{R}\cup \{\infty \}$ is identified with the boundary $\partial{\mathbb{H}^2}$ of $\mathbb{H}^2$.   

Observe that in both the cases it is possible to determine the axis of $M_\gamma$ (or $M_\delta$) from the points where they intersect $\mathbb{R}$. Now these last set of points are just the fixed points of $M_\gamma$ (or $M_\delta$). The fixed points of $M_\gamma$ can be computed simply as follows. There are two cases.

\textbf{Case I:} First let the axis of $M_\gamma$ (or $M_\delta$) be a semi-circle. Then both the points of intersections are finite real numbers that satisfy
\begin{equation*}
M_\gamma(z) = z \Rightarrow \frac{a_\gamma z + b_\gamma}{c_\gamma z + d_\gamma} = z \Rightarrow c_\gamma z^2 + (d_\gamma - a_\gamma)z - b_\gamma = 0.
\end{equation*}
Hence the two points of intersections of the axis of $M_\gamma$ with the real line are the two roots of the equation
\begin{equation}
c_\gamma z^2 + (d_\gamma - a_\gamma)z - b_\gamma = 0.
\end{equation}
Denote these by $\alpha_\gamma$ and $\beta_\gamma$ with $\alpha_\gamma < \beta_\gamma$. In terms of matrix coefficients of $M_\gamma$ we have
\begin{equation}
\alpha_\gamma = \frac{(a_\gamma - d_\gamma) - \sqrt{(a_\gamma - d_\gamma)^2 + 4c_\gamma b_\gamma}}{2 c_\gamma}, $$$$ \beta_\gamma = \frac{(a_\gamma - d_\gamma) + \sqrt{(a_\gamma - d_\gamma)^2 + 4c_\gamma b_\gamma}}{2 c_\gamma}.
\end{equation}
Using $\det{M_\gamma} = 1$ they take the form:
\begin{equation}
\alpha_\gamma = \frac{(a_\gamma - d_\gamma) - \sqrt{\tr{M_\gamma}^2 - 4}}{2 c_\gamma}, ~~ \beta_\gamma = \frac{(a_\gamma - d_\gamma) + \sqrt{\tr{M_\gamma}^2 - 4}}{2 c_\gamma}.
\end{equation}
In particular the center and the Euclidean radius of the axis of $M_\gamma$ are respectively
\begin{equation}
m_\gamma = (\frac{a_\gamma - d_\gamma}{2c_\gamma}, 0) ~~ ~~ \textrm{and} ~~ ~~ r_\gamma = \frac{\sqrt{\tr{M_\gamma}^2 - 4}}{2c_\gamma}.
\end{equation}
\textbf{Case II:} The axis of $M_\gamma$ is a vertical straight line. In particular $c_\gamma = 0$. Then the only point of intersection between the axis of $M_\gamma$ and $\mathbb{R}$ is $(\frac{b_\gamma}{d_\gamma - a_\gamma}, 0)$.
\subsection{Cosine of the angle}
Consider two intersecting closed geodesics $\gamma$ and $\delta$ on $S$. Fix one point of their intersection $p$. Choose two representatives $M_\gamma, M_\delta$ for $\gamma, \delta$ respectively such that the point of intersection $\tilde{p}$ between the axis of $M_\gamma$ and the axis of $M_\delta$ is a lift of $p$ under the covering map $\pi: \mathbb{H}^2 \to \mathbb{H}^2/\Gamma =S$. Let $\theta = \theta(\gamma, \delta, p)$. Hence $\theta$ is the angle between the axis of $M_\gamma$ and the axis of $M_\delta$ at $\tilde{p}$. Now we have two cases depending on the nature of the axes of $M_\gamma$ and $M_\delta$. We treat them separately.

\begin{figure}[H]
	\includegraphics[scale=.45]{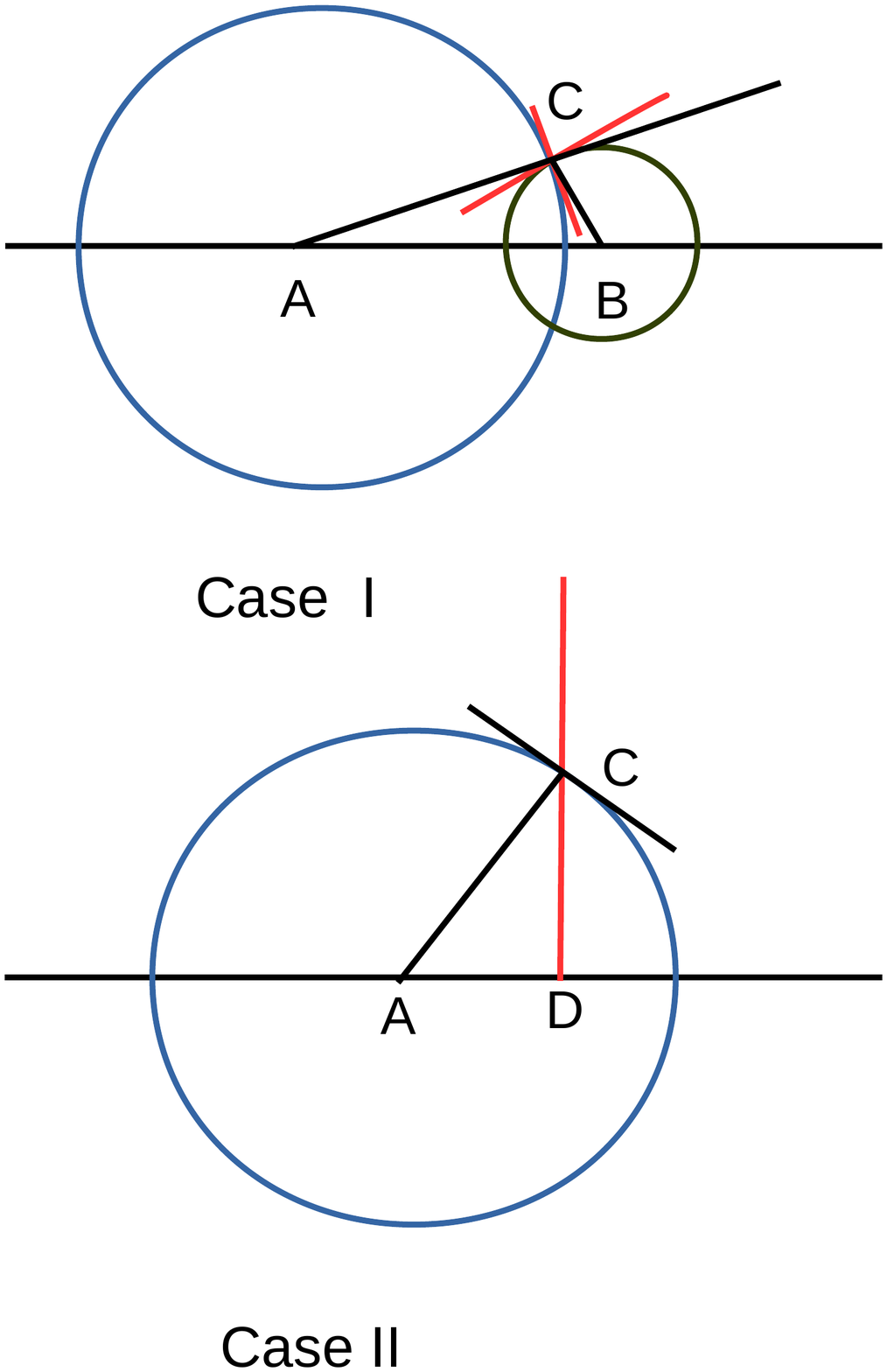}
	\caption{Cosine of the angle}
	\label{fig:angle-cosine}
\end{figure}

\textbf{Case I:}  First let us assume that both $M_\gamma$ and $M_\delta$ have semi-circle axes. This situation is explained in the top picture in \ref{fig:angle-cosine}. Let $\psi$ be the angle between the normals to the the axis of $M_\gamma$ and the axis of $M_\delta$ at $\tilde{p}$. Then $\psi = \pi - \theta$. 

Now consider the Euclidean triangle on $\mathbb{H}^2$ with the following three vertices: the centre of the (semi-circle) axis of $M_\gamma$, the centre of the (semi-circle) axis of $M_\delta$ and $\tilde{p}$ the point of intersection of the two axes. Let us denote the distance between the two centres by $d_{\gamma,\delta}$. Hence
\begin{equation}
d_{\gamma, \delta} = |\frac{a_\gamma - d_\gamma}{2c_\gamma} - \frac{a_\delta - d_\delta}{2c_\delta}|.
\end{equation}
Using Euclidean geometry for the above described triangle we obtain
\begin{equation*}
\cos(\pi-\theta) = \frac{r_\gamma^2 + r_\delta^2 - d_{\gamma, \delta}^2}{2r_\gamma r_\delta}.
\end{equation*}
Thus
\begin{equation}\label{eucld1}
\cos^2(\theta) = \frac{(r_\gamma^2 + r_\delta^2 - d_{\gamma, \delta}^2)^2}{4r_\gamma^2 r_\delta^2}.
\end{equation}
\textbf{Case II:} Now we assume that the axis of $M_\gamma$ is a vertical straight line. Since $\gamma$ and $\delta$ intersect each other, the axis of $M_\delta$ must be a semi-circle. This situation is explained in the bottom picture of \ref{fig:angle-cosine}. Consider the normal $N_\delta$ to the axis of $M_\delta$ at $\tilde{p}$. Let $\psi$ be the angle between the boundary $\partial{\mathbb{H}^2} = \mathbb{R}$ and $N_\delta$. Observe that $\psi = \theta$. Now we consider the Euclidean triangle with vertices: the center of the axis of $M_\delta$, the point of intersection between axis of $M_\gamma$ and $\partial{\mathbb{H}^2} = \mathbb{R}$ and $\tilde{p}$. By the definition of the cosine function and the last equality $\psi = \theta$ we get
\begin{equation}\label{eucld2}
	\cos^2(\theta) = \bigg(\frac{|\frac{b_\gamma}{d_\gamma - a_\gamma} - \frac{a_\delta - d_\delta}{2c_\delta}|}{\frac{\sqrt{\tr{M_\delta}^2 - 4}}{2c_\delta}}\bigg)^2.
\end{equation}

\begin{rem}
	From the above two expressions it is clear that $\cos^2(\theta)$ is expressible as rational functions in the matrix entries $a_\gamma, b_\gamma, c_\gamma, d_\gamma$ of $M_\gamma$ and $a_\delta, b_\delta, c_\delta, d_\delta$ of $M_\delta$.
\end{rem} 
\section{Proof of Theorem \ref{fuchsian}}
Consider the field 
$$\mathbb{F}_\Gamma = \mathbb{Q}(a_\gamma, b_\gamma, c_\gamma, d_\gamma: M_\gamma = \left( \begin{array}{cc}
a_\gamma & b_\gamma \\
c_\gamma & d_\gamma \
\end{array} \right) \in \Gamma)$$ 
generated by the entries of the matrices in $\Gamma \subset \textrm{SL}(2, \mathbb{R})$. Observe that this is a finitely generated field. This is clear because $\Gamma$ is a finitely generated group and so adjoining the matrix entries of a generating subset of $\Gamma$ is enough. 

Now we have two cases. The filed $\mathbb{F}_\Gamma$ may or may not be algebraic over $\mathbb{Q}$. Let us assume that $\mathbb{F}_\Gamma$ is not algebraic over $\mathbb{Q}$. Since $\mathbb{F}_\Gamma$ is finitely generated over $\mathbb{Q}$, there is a purely transcendental extension ${\mathbb{T}_\Gamma}{|_\mathbb{Q}} \subset \mathbb{F}_\Gamma$ such that $\mathbb{F}_\Gamma$ is algebraic over $\mathbb{T}_\Gamma$. To treat the two cases at the same time let ${\mathbb{T}_\Gamma}$ denote $\mathbb{Q}$ when $\mathbb{F}_\Gamma$ is algebraic over $\mathbb{Q}$. In both the cases $\mathbb{F}_\Gamma$ is finitely generated over $\mathbb{T}_\Gamma$. Hence the degree $[\mathbb{F}_\Gamma: \mathbb{T}_\Gamma]$ of the extension ${\mathbb{F}_\Gamma}{|_{\mathbb{T}_\Gamma}}$ is finite.

Now from \eqref{eucld1} and \eqref{eucld2} it is clear that for $\theta = \theta(\gamma, \delta, p)$ the values $\cos^2(\theta) \in \mathbb{F}_\Gamma$. Hence the degree of the field extension $\mathbb{F}_\Gamma(e^{2i\theta})_{|{\mathbb{F}_\Gamma}}$
\begin{equation*}
 [\mathbb{F}_\Gamma(e^{2i\theta}):{\mathbb{F}_\Gamma}] \le 2.
\end{equation*}
This implies that the degree of the extension $\mathbb{T}_\Gamma(e^{2i\theta})_{|{\mathbb{T}_\Gamma}}$
\begin{equation*}
 [\mathbb{T}_\Gamma(e^{2i\theta}): \mathbb{T}_\Gamma] \le [\mathbb{F}_\Gamma(e^{2i\theta}): \mathbb{T}_\Gamma]  = [\mathbb{F}_\Gamma(e^{2i\theta}):{\mathbb{F}_\Gamma}]\cdot[\mathbb{F}_\Gamma: \mathbb{T}_\Gamma] \le 2[\mathbb{F}_\Gamma: \mathbb{T}_\Gamma].
\end{equation*}

Now recall that $\mathbb{T}_\Gamma$ is a purely transcendental extension of $\mathbb{Q}$ and so for $\theta$ rational multiple of $\pi$ (since $e^{2i\theta}$ is algebraic over $\mathbb{Q}$) we always have $$[\mathbb{T}_\Gamma(e^{2i\theta}): \mathbb{T}_\Gamma] = [\mathbb{Q}(e^{2i\theta}): \mathbb{Q}].$$ Now let $\theta = \frac{p}{q}\pi$. It is a know fact that the degree $[\mathbb{Q}(e^{2i\theta}): \mathbb{Q}] = \phi(q)$ where $\phi$ is the Euler $\phi$-function. Thus combining the above inequalities we have $$\phi(q) \le 2[\mathbb{F}_\Gamma: \mathbb{T}_\Gamma].$$ Hence there are only finitely many choices for $q$ by \cite[Theorem 328]{H-W}.
\begin{rems}
	(i) Observe that the field $\mathbb{F}_\Gamma$ depends explicitly on $\Gamma$ where $A(\mathbb{H}^2/\Gamma)$ depends only on the conjugacy class of $\Gamma$ because conjugate groups produce isometric surfaces. Hence we conclude that for any $\frac{p}{q}\cdot \pi \in A(\mathbb{H}^2 /\Gamma)$ $$ \phi(q) \le 2 \cdot \min_{\gamma \in \textrm{PSL}(2, \mathbb{R})} [\mathbb{F}_{\gamma \Gamma {\gamma^{-1}}}, \mathbb{T}_{\gamma \Gamma {\gamma^{-1}}}].$$ This can be used to give an explicit bound on the size of $A(\mathbb{H}^2 /\Gamma) \cap \mathbb{Q}\cdot \pi.$
	
	(ii) For the modular surface $\mathbb{H}^2/ \text{PSL}(2, \mathbb{Z})$ the group $\Gamma$ is $\text{PSL}(2, \mathbb{Z})$ and so the field $\mathbb{F}_\Gamma$ is just $\mathbb{Q}$. Hence for any $\frac{p}{q} \pi \in A(\mathbb{H}^2/\text{PSL}(2, \mathbb{Z}))$ we have $\phi(q) \le 2$ i.e. $q \le 6.$ A simple computation provides that the possible angles are $\pi/6, \pi/4$ and $\pi/3$. %It is not very difficult to see that there actually are closed geodesics on $\mathbb{H}^2/\text{PSL}(2, \mathbb{Z})$ that intersect each other with those angles.
\end{rems}
\section{Some questions and examples}
 Let $\Gamma$ be a Fuchsian group as above and $S = \mathbb{H}^2/\Gamma$. Given an angle $\theta \in A(S)$ one may consider the map $\Theta: A(S) \to \mathbb{F}_\Gamma^1$ given by $\Theta(\theta) = \cos^2(\theta)$, where $\mathbb{F}_\Gamma^1 \subset \mathbb{F}_\Gamma$ is the set of elements with norm $< 1$. 
 \begin{que}
 	What is the image of this map ? 
 \end{que}
  Since $A(S)$ is dense in $[0, \pi]$ the image is dense in $[-1, 1]$ and hence in $\mathbb{F}_\Gamma^1$. It is not clear if it equals  $\mathbb{F}_\Gamma^1$ though.
\subsection{An angle with infinite multiplicity}  
Now we show that there are closed hyperbolic surfaces $S$ such that $\pi/2$ has infinite multiplicity in $\mathcal{A}(S)$. To construct such a surface start with a compact surface $S'$ with geodesic boundary and consider its double $DS'$. Observe that $DS'$ has a reflection isometry along the boundary geodesics of $S'$. Now take a closed geodesic that is symmetric with respect to this reflection and intersects at least one of the boundary geodesics of $S'$. By the reflection symmetry each of these angles of intersection must be equal to $\pi/2$. It is not hard to construct a surface having infinitely many geodesics of this type. This poses our last question.
\begin{que}
What angles in $\mathcal{A}(S)$ can have infinite multiplicity ? 
\end{que}
%It seems reasonable to believe that if there is an angle with infinite multiplicity then the surface has an isometry. This is our last question.
 %\begin{que}
 %	Does the existence of such an angle of infinite multiplicity imply that the surface has a symmetry ?
 %\end{que}

\end{document}